\documentclass[10pt]{article}

\usepackage{amscd,amsmath, amssymb, fancyhdr, mathbbol}

\numberwithin{equation}{section}

\newcommand{\version}{version 3.0,\ \ July 18, 2016}

\def\eqref#1{(\ref{#1})}
\newcommand{\goth}{\mathfrak}

\newcommand{\arrow}{{\:\longrightarrow\:}}
\newcommand{\Z}{{\Bbb Z}}
\newcommand{\C}{{\Bbb C}}
\newcommand{\R}{{\Bbb R}}
\newcommand{\Q}{{\Bbb Q}}

\def\1{\sqrt{-1}\:}
\newcommand{\restrict}[1]{{\left|_{{\phantom{|}\!\!}_{#1}}\right.}}
\newcommand{\cntrct}                
{\hspace{2pt}\raisebox{1pt}{\text{$\lrcorner$}}\hspace{2pt}}

\newcommand{\calo}{{\cal O}}

\renewcommand{\bar}{\overline}
\renewcommand{\phi}{\varphi}
\renewcommand{\epsilon}{\varepsilon}
\renewcommand{\geq}{\geqslant}
\renewcommand{\leq}{\leqslant}

\newcommand{\Teich}{\operatorname{Teich}}

\newcommand{\diam}{\operatorname{\sf diam}}

\newcommand{\Vol}{\operatorname{Vol}}
\newcommand{\Pic}{\operatorname{Pic}}

\newcommand{\Sym}{\operatorname{Sym}}

\newcommand{\Aut}{\operatorname{Aut}}

\newcommand{\Diff}{\operatorname{Diff}}
\newcommand{\Spin}{\operatorname{Spin}}

\newcommand{\Tr}{\operatorname{Tr}}

\newcommand{\Comp}{\operatorname{Comp}}

\newcounter{Mycounter}[section]
\newcounter{lemma}[section]
\renewcommand{\thelemma}{{Lemma \thesection.\arabic{lemma}}}
\newcommand{\lemma}{%
   \setcounter{lemma}{\value{Mycounter}}
   \refstepcounter{lemma}
   \stepcounter{Mycounter}
   {\noindent \bf \thelemma:\ }}

\newcounter{claim}[section]
\renewcommand{\theclaim}{{Claim \thesection.\arabic{claim}}}
\newcommand{\claim}{%
   \setcounter{claim}{\value{Mycounter}}
   \refstepcounter{claim}
   \stepcounter{Mycounter}
   {\noindent \bf \theclaim:\ }}

\newcounter{sublemma}[section]
\newcounter{corollary}[section]
\renewcommand{\thecorollary}{{Corollary \thesection.\arabic{corollary}}}
\newcommand{\corollary}{%
   \setcounter{corollary}{\value{Mycounter}}
   \refstepcounter{corollary}
   \stepcounter{Mycounter}
   {\noindent \bf \thecorollary:\ }}

\newcounter{theorem}[section]
\renewcommand{\thetheorem}{{Theorem \thesection.\arabic{theorem}}}
\newcommand{\theorem}{%
   \setcounter{theorem}{\value{Mycounter}}
   \refstepcounter{theorem}
   \stepcounter{Mycounter}
   {\noindent \bf \thetheorem:\ }}

\newcounter{conjecture}[section]
\renewcommand{\theconjecture}{{Conjecture 
\thesection.\arabic{conjecture}}}
\newcommand{\conjecture}{%
   \setcounter{conjecture}{\value{Mycounter}}
   \refstepcounter{conjecture}
   \stepcounter{Mycounter}
   {\noindent \bf \theconjecture:\ }}

\newcounter{proposition}[section]
\renewcommand{\theproposition}
     {{Proposition \thesection.\arabic{proposition}}}
\newcommand{\proposition}{%
   \setcounter{proposition}{\value{Mycounter}}
   \refstepcounter{proposition}
   \stepcounter{Mycounter}
   {\noindent \bf \theproposition:\ }}

\newcounter{definition}[section]
\renewcommand{\thedefinition}
     {{Definition~\thesection.\arabic{definition}}}
\newcommand{\definition}{%
   \setcounter{definition}{\value{Mycounter}}
   \refstepcounter{definition}
   \stepcounter{Mycounter}
   {\noindent \bf \thedefinition:\ }}

\newcounter{example}[section]
\newcounter{remark}[section]
\renewcommand{\theremark}{{Remark \thesection.\arabic{remark}}}
\newcommand{\remark}{%
   \setcounter{remark}{\value{Mycounter}}
   \refstepcounter{remark}
   \stepcounter{Mycounter}
   {\noindent \bf \theremark:\ }}

\newcounter{problem}[section]
\newcounter{question}[section]
\makeatletter

\@addtoreset{equation}{section} \@addtoreset{footnote}{section}
\makeatother

\def\blacksquare{\hbox{\vrule width 5pt height 5pt depth 0pt}}
\def\endproof{\blacksquare}

\begin{document}
\begin{center}
{\LARGE\bf
On the Kobayashi pseudometric, complex automorphisms and 
hyperk\"ahler manifolds\\[4mm]
}

Fedor Bogomolov\footnote{Partially supported by a Simons Travel grant}, 
Ljudmila Kamenova\footnote{Partially supported by NSF DMS-1502154}, 
Steven Lu\footnote{Partially supported by an NSERC discovery grant}, 
Misha Verbitsky\footnote{Partially supported 
by RScF grant, project 14-21-00053, 11.08.14.}

\end{center}

{\small \hspace{0.10\linewidth}
\begin{minipage}[t]{0.85\linewidth}
{\bf Abstract} \\
We define the Kobayashi quotient of a complex variety by identifying 
points with vanishing Kobayashi pseudodistance between them and show that if a complex projective manifold has an automorphism whose order is infinite, then the fibers of this quotient map are nontrivial. 
We prove that the Kobayashi quotients associated to ergodic complex 
structures on a compact manifold are isomorphic. 
We also give a proof of Kobayashi's conjecture 
on the vanishing of the pseudodistance 
for hyperk\"ahler manifolds having Lagrangian fibrations without multiple fibers 
in codimension one. For a hyperbolic automorphism of a hyperk\"ahler manifold, we prove that its cohomology eigenvalues are determined by its Hodge numbers, compute its dynamical 
degree and show that its cohomological trace grows exponentially, giving estimates on the number of its periodic points. 
\end{minipage}
}

\tableofcontents


\section{Introduction}


Kobayashi conjectured that a compact K\"ahler manifold with 
semipositive Ricci curvature has vanishing Kobayashi pseudometric. 
In a previous paper (\cite{klv}) Kamenova-Lu-Verbitsky have proved the 
conjecture for all K3 surfaces and for certain hyperk\"ahler manifolds 
that are deformation equivalent to Lagrangian fibrations. Here we give 
an alternative proof of this conjecture for hyperk\"ahler 
Lagrangian fibrations without multiple fibers in codimension one, see Section~\ref{S_Abelian}. 

\hfill

\theorem
Let $f\colon M \arrow B = \C{\Bbb P}^n$ be a hyperk\"ahler Lagrangian 
fibration 
without multiple fibers in codimension one over $B$. Then the Kobayashi pseudometric $d_M$ vanishes identically on $M$ and the Royden-Kobayashi pseudonorm $|\ |_M$ vanishes identically on a Zariski open subset of $M$. 

\hfill

In Section~\ref{S_Infinite}, we explore compact complex manifolds $M$ having an automorphism of infinite order. If such a manifold is projective, we show that the Kobayashi pseudometric is 
everywhere degenerate. For each point $x \in M$ we define the subset 
$M_x \subset M$ of points in $M$ whose pseudo-distance to $x$ is zero. 
Define the relation $x\sim y$ on $M$ given by $d_M(x,y)= 0$. 
There is a well defined set-theoretic quotient map 
$\Psi : M \arrow S= M/\!\sim$, called {\bf the
Kobayashi quotient map}. We say that $|\ |_M$ is {\bf Voisin-degenerate} 
at a point $x\in M$ if 
there is a sequence of holomorphic maps $\phi_n: D_{r_n}\to M$ such that 
$\phi_n(0)\to x,\ |\phi_n^{\prime}(0)|_h=1\ \mbox{and}\ r_n\to \infty.$

\hfill


\theorem
Let $M$ be a  
complex projective 
manifold with an automorphism $f$
of infinite order. Then the Kobayashi pseudo-metric $d_M$ 
is everywhere degenerate in the sense that $M_x\neq \{x\}$ for
all $x\in M$. The Royden-Kobayashi pseudo-norm $|\ |_M$ 
is everywhere Voisin-degenerate. Moreover, 
every fiber of the map $\Psi: M\arrow S$ constructed above 
contains a Brody curve and is connected. 

\hfill

Define {\bf the Kobayashi quotient} $M_K$ of $M$ to be 
the space of all equivalence classes 
$\{x\sim y \ \ |\ \ d_M(x,y)=0\}$ equipped with 
the metric induced from $d_M$. 

In Section~\ref{S_Ergoquotients}, we show that the Kobayashi quotients for 
ergodic complex structures are isometric, equipped with  
the natural quotient pseudometric. This generalizes the 
key technical result of \cite{klv} for the identical vanishing of $d_M$
for ergodic complex structures on hyperk\"ahler manifolds.

\hfill

\theorem
Let $(M,I)$ be a compact complex manifold, and $(M,J)$ its deformation.
Assume that the complex structures $I$ and $J$ are both ergodic.
Then the corresponding Kobayashi quotients are isometric.

\hfill

Finally in Section~\ref{S_Eigen}, we prove that the cohomology eigenvalues of 
a hyperbolic automorphism of a hyperk\"ahler manifold 
are determined by its Hodge numbers. We compute its dynamical 
degree in the even cases and give an upper bound in the odd cases. 

\hfill

\theorem
Let $(M,I)$ be a hyperk\"ahler manifold, and
$T$ a hyperbolic automorphism acting on cohomology
as $\gamma$. Denote by $\alpha$ the
eigenvalue of $\gamma$ on $H^2(M,\R)$ with $|\alpha|>1$.
Then all eigenvalues of
$\gamma$ have absolute value which is a power of
$\alpha^{1/2}$. Moreover, the maximal of these eigenvalues
on even cohomology $H^{2d}(M)$ is equal to $\alpha^{d}$, and finally, 
on odd cohomology $H^{2d+1}(M)$ the maximal eigenvalue of $\gamma$ is 
strictly less than $\alpha^{\frac{2d+1}{2}}$.

\hfill

As a corollary we obtain that the trace $\Tr(\gamma^N)$ grows 
asymptotically as $\alpha^{nN}$. We also show that the number of 
$k$-periodic points grows as $\alpha^{nk}$.

\hfill

The work on this paper started during the Simons Symposium 
``Geometry over nonclosed fields'' held in March, 2015. 
The authors are grateful to the Simons Foundation for  
providing excellent research conditions.


\section{Preliminaries}


\definition
A {\em hyperk\"ahler} (or {\em irreducible holomorphic symplectic}) 
manifold 
$M$ is a compact complex K\"ahler manifold with $\pi_1(M)=0$ and 
$H^{2,0}(M)=\C \sigma$ where $\sigma$ is everywhere non-degenerate. 

\hfill

Recall that a fibration
is a connected surjective holomorphic map.
On a hyperk\"ahler manifold the structure of a fibration, if one
exists, is limited by Matsushita's theorem. 

\hfill

\theorem (Matsushita, \cite{_Matsushita:fibred_})
Let $M$ be a hyperk\"ahler manifold and $f\colon M\arrow B$ a fibration 
with 
$0 < \dim B < \dim M$.
Then $\dim B =\frac{1}{2}\dim M$ and the general fiber of $f$ is a 
Lagrangian abelian variety. The base $B$ has at worst $\Q$-factorial 
log-terminal singularities, has Picard number $\rho(B)=1$ and $-K_B$ is 
ample. 

\hfill

\remark
$B$ is smooth in all of the known examples. 
It is conjectured that $B$ is 
always smooth. 

\hfill

\theorem(Hwang \cite{Hwang})
In the settings above, if $B$ is smooth then $B$ is isomorphic to 
$\C{\Bbb P}^n$, where $\dim_{\C} M = 2n$.

\hfill

\definition
Given a hyperk\"ahler manifold $M$, there is a non-degenerate 
integral quadratic form $q$ on 
$H^2(M,\Z)$, called the {\em Beauville-Bogomolov-Fujiki form} 
(BBK form for short), 
of signature $(3,b_2-3)$ and satisfying the {\em Fujiki relation} 
$$\int_M \alpha^{2n}=c\cdot q(\alpha)^n\qquad\text{for }\alpha \in 
H^2(M,\Z),$$ 
with $c>0$ a constant depending on the topological type of $M$. 
This form generalizes the intersection pairing on K3 surfaces. 
For a detailed description of the form we refer the reader to 
\cite{_Fujiki:HK_},  \cite{_Beauville_} and \cite{_Bogomolov:defo_}. 

\hfill

\remark
Given $f\colon M \arrow \C{\Bbb P}^n$, $h$ the hyperplane 
class on $\C \Bbb P^n$, and $\alpha=f^*h$, 
then $\alpha$ is nef and $q(\alpha)=0$.

\hfill

\conjecture[SYZ]
If $L$ is a nontrivial nef line bundle on $M$ with $q(L)=0$, 
then $L$ induces a Lagrangian fibration, given as above.

\hfill

\remark
This conjecture is known for deformations of Hilbert schemes of points on 
K3
surfaces 
(Bayer--Macr\`i \cite{_Bayer_Macri_}; Markman \cite{_Markman:SYZ_}), 
and for deformations of the generalized Kummer varieties $K_n(A)$ 
(Yoshioka \cite{_Yoshioka_}). 


\hfill

\definition 
The {\em Kobayashi pseudometric} on $M$ is the maximal 
pseudometric $d_M$ such that all holomorphic maps 
$f\colon (D,\rho) \arrow (M,d_M)$ are distance 
decreasing, where $(D,\rho)$ is the unit disk with the Poincar\'e metric. 

\hfill

\definition
A manifold $M$ is {\em Kobayashi hyperbolic} if $d_M$ is a metric, 
otherwise 
it is called {\em Kobayashi non-hyperbolic}. 

\hfill



\remark In \cite{_Kobayashi:1976_}, it is asked whether a compact K\"ahler 
manifold $M$ of semipositive Ricci curvature has identically vanishing 
pseudometric, which we denote by $d_M\equiv 0$. The question applies to 
hyperk\"ahler manifolds but was unknown even for the case of  surfaces 
outside the projective case. But Kamenova-Lu-Verbitsky (in \cite{klv}) 
have recently resolved completely the case of surfaces with the following affirmative results. \\[-1mm]

\theorem \cite{klv}
Let $S$ be a K3 surface. Then $d_S\equiv 0$.\\[-1mm]

\remark  A birational version of a conjecture of Kobayashi \cite{_Kobayashi:1976_} would state that a compact hyperbolic manifold be of general type if its Kobayashi pseudometric is nondegenerate somewhere (i.e. nondegenerate on some open set). This was open for surfaces but now resolved outside surfaces of class VII.

\hfill

\theorem \cite{klv} \label{_klv_}
Let $M$ be a hyperk\"ahler manifold of non-maximal Picard rank 
and deformation 
equivalent to a Lagrangian fibration. Then $d_M\equiv 0$.

\hfill

\theorem \cite{klv} \label{klvdouble}
Let $M$ be a hyperk\"ahler manifold with $b_2(M)\ge 7$ (expected to always hold) and 
with 
maximal Picard rank $\rho = b_2-2$. 
Assume the SYZ conjecture for deformations of $M$. Then $d_M\equiv 0$.

\hfill

\remark Except for the proof of \ref{klvdouble},
we indicate briefly a proof of these theorems below. \ref{klvdouble}
is proved in \cite{klv} using the existence of double Lagrangian 
fibrations
on certain deformations of $M$. 
Here we give a different proof of vanishing of the Kobayashi 
pseudometric for certain hyperk\"ahler Lagrangian fibrations without 
using double fibrations. 

\hfill

\definition
Let $M$ be a compact complex manifold and 
$\Diff^0(M)$ the connected component to identity of its diffeomorphism 
group. 
Denote by $\Comp$ the space of complex structures on $M$, equipped with 
a structure of Fr\'echet manifold. The {\em Teichm\"uller space} of $M$ 
is the quotient  $\Teich:=\Comp/\Diff^0(M)$. 
The Teichm\"uller space is finite-dimensional for $M$ Calabi-Yau 
(\cite{_Catanese:moduli_}). 
Let $\Diff^+(M)$ be the group of orientable diffeomorphisms of 
a complex manifold $M$. The {\em mapping class group} 
$\Gamma:=\Diff^+(M)/\Diff^0(M)$ acts on $\Teich$. 
An element $I \in \Teich$ is called {\em ergodic} if the orbit 
$\Gamma \cdot I$ is dense in $\Teich$, where 
\[
\Gamma \cdot I = \{I' \in \Teich : (M,I)\sim (M,I')\}.
\]

\hfill

\theorem(Verbitsky, \cite{_Verbitsky:ergodic_})
If $M$ is hyperk\"ahler and $I\in \Teich$, then $I$ is ergodic 
if and only if $\rho(M,I)<b_2-2.$

\hfill

\remark For a K3 surface $(M,I)$ not satisfying the above condition on the
Picard rank $\rho$, it is easily seen to admit Lagrangian (elliptic) 
fibrations 
over $\C{\Bbb P}^1$ without multiple fibers, and it is projective. 
Then $d_{(M,J)}\equiv 0$ by \ref{nomultiple} below, for example.

\hfill

\proposition\label{usc}
Let $(M,J)$ be a compact complex manifold with $d_{(M,J)}\equiv 0$.
Let $I \in \Teich$ be an ergodic complex structure deformation equivalent 
to 
$J$. Then $d_{(M,I)} \equiv 0$. 

\hfill

{\bf Proof:} Here we shall reproduce the proof from \cite{klv}. 
Consider the diameter function 
$\diam : \; \Teich \arrow \R_{\geq 0},$ 
the maximal distance between two points. It is upper semi-continuous 
(Corollary 1.23 in \cite{klv}). Since the complex structure $J$ is in the 
limit set 
of the orbit of the ergodic structure $I$, by upper semi-continuity 
$0 \leq \diam(I) \leq \diam(J) =0$. \endproof


\section{(Royden-)Kobayashi pseudometric on Abelian fibrations}\label{S_Abelian}



The following lemma is a generalization of Lemma 3.8 in 
\cite{_Bu_Lu_} to the case of abelian fibrations. The generalization is given for example in the 
Appendix of \cite{klv}. Recall that an abelian fibration is a connected 
locally projective surjective K\"ahler morphism with abelian varieties 
as fibers. 

\hfill

\lemma \label{_BL_gen_}
Let $\pi\colon T \arrow C$ be an abelian fibration over a non-compact 
complex curve $C$ which locally has sections and such that not all 
components of the fibers are multiple. Then $T$ has an analytic section 
over $C$. This is the case if $\pi$ has no multiple fibers.

\hfill

{\bf Proof:} 
There is a Neron model $N$ for $T$ and a short exact sequence 
$$0 \arrow F \arrow \calo(L) \arrow \calo(N) \arrow 0$$ 
where $L$ is a vector bundle, $F$ is a sheaf of groups $\Z^{2n}$ with 
degenerations, i.e., sheaf of discrete subgroups with generically maximal 
rank, and $\calo(N)$ is the sheaf of local sections of $N$ 
(whose general fibers are abelian varieties). 
Thus $T$ corresponds to an element $\theta$ in  $H^1(C,\calo(N))$. There is 
an 
induced exact sequence of cohomologies: 
$H^1(C,\calo(L))\arrow H^1(C,\calo(N))\arrow H^2(C,F)$. 
Note that $H^1(C,\calo(L))= 0$ since $C$ is Stein, and $H^2(C,F)= 0$ 
since it is topologically one-dimensional. 
Thus $\theta=0$ and hence there is an analytic section. 
The last part of the lemma is given by Proposition~4.1 of \cite{klv}. 
\endproof 

\hfill


\theorem\label{nomultiple}
Let $f\colon M \arrow B = \C{\Bbb P}^n$ be a hyperk\"ahler Lagrangian 
fibration 
without multiple fibers in codimension one over $B$. Then $d_M\equiv 0$ and $|\ |_M$ vanishes on a nonempty Zariski open subset of $M$. 

\hfill

{\bf Proof:} 
The fibers of $f$ are projective, and furthermore, there is a canonical 
polarization on them (see \cite{_Oguiso_ST_} and \cite{_Oguiso_rank_1_}, 
respectively). 
This also follows from 
\cite{_Verbitsky:degenerate_twistor_}, Theorem 1.10, which
implies that the given fibration is diffeomorphic to another
fibration $f\colon M' \arrow B$ with holomorphically the same fibers and 
the same base, but with projective total space $M'$. 
Standard argument (via the integral lattice in the ``local'' 
Neron-Severi group) now shows that $f$ is locally projective. 

By assumption, there are no multiple fibers outside a codimension~$2$ 
subset $S \subset B$ whose complement $U$ contains at most the smooth
codimension-one part $D_0$ of the discriminant locus of $f$ where 
multiplicity
of fibers are defined locally generically. Since the pseudometric is  
unchanged after removing codimension $2$ subsets 
(\cite{_Kobayashi:1998_}), 
it is enough to restrict the fibration to that over $U$. 

Let $C={\Bbb P}^1$ be a line in $B= {\Bbb P}^n$ contained in $U$
(and intersecting $D_0$ transversely). Then $f$ restricts to an abelian 
fibration
$X=f^{-1}(C)$ over $C$ without multiple fibers and so \ref{_BL_gen_}
applies to give a section over the affine line 
$A^1 = C\setminus(\infty)$. 

As $S$ is codimension two or higher, we can connect any 
two general points in $U$ by a chain of such 
$A^1$'s in $U$. One can thus connect two general 
points $x$ and $y$ on $M$ by a chain consisting of fibers and sections 
over 
the above $A^1$'s. Since the Kobayashi 
pseudometric vanishes on each fiber and each such section, the triangle 
inequality 
implies $d_M(x,y)=0$. Therefore $d_M$ vanishes 
on a dense open subset of $M$ and hence
$d_M\equiv 0$ by the continuity of $d_M$. 

The same argument gives the vanishing statement of $|\ |_M$ 
 via Theorem A.2 of \cite{klv}.
\endproof

\hfill

\remark
In the theorem above, it is sufficient to assume that $B$ is nonsingular and that $d_B\equiv 0$, true if $B$ is rationally connected. In fact, if one assumes further the vanishing of 
$|\ |_B$ on a nonempty Zariski open, then the 
same is true for $|\ |_M$,
generalizing the corresponding theorems in \cite{klv}. 
The reader should have no difficulty to see these by
the obvious modifications of the above proof.

\section{Automorphisms of infinite order}\label{S_Infinite}


We first sketch the proof of Kobayashi's theorem that Kobayashi hyperbolic 
manifolds have only finite order automorphisms 
(Theorem 9.5 in \cite{_Kobayashi:1976_}). 

\hfill

\theorem
Let $M$ be a Kobayashi hyperbolic manifold.
Then its group of birational transformations is finite.

\hfill

{\bf Proof:} 
First, notice that a birational self-map is a composition of
a blow-up, an automorphism and a blow-down. Since $M$ contains
no rational curves, any birational self-map is holomorphic,
and we need to prove the finiteness of the automorphism group.

Observe that the automorphisms of 
a hyperbolic manifold are isometries of the Kobayashi metric. Also 
the group of isometries of a compact metric space is compact with respect
to the compact open topology by a theorem of 
Dantzig and Van der Waerden, see for example 
\cite[Theorem~5.4.1]{_Kobayashi:1998_}. 
On the other hand, compact Kobayashi 
hyperbolic manifolds have no holomorphic vector fields, because 
each such vector field gives an orbit which is an entire curve. This 
means that the group of holomorphic automorphisms $\text{Aut}(M)$
 of $M$ is discrete as it is a complex Lie group in the compact open topology acting holomorphically on $M$ by the work of Bochner-Montgomery
 \cite{BM1,BM2}. Since $\text{Aut}(M)$ 
is discrete and compact, this means it is finite. 
\endproof

\hfill

Consider the pseudo-distance function $d_M: M \times M \arrow {\Bbb R}$, 
defined by the Kobayashi pseudo-distance $d_M(x,y)$ on pairs $(x,y)$. It 
is 
a symmetric continuous function which is bounded for compact $M$. 
Since it is symmetric, we can consider $d_M$ as a function on the 
symmetric product $\Sym^2 M$ with $d_M= 0$ on the diagonal. 

\hfill

\lemma 
There is a compact space $S$ with a continuous map 
$\Psi : M \arrow S$ and there is a distance function 
$d_S$ on $S$ making $S$ into a 
compact metric space 
such that $d_M = d_S \circ \psi$,
where $\psi : \Sym^2 M \arrow \Sym^2 S$ 
is the map induced by $\Psi$.

\hfill

{\bf Proof:} The subset $ M_x \subset M$ of points 
$y \in M$ with $d_M(x,y)= 0$ is compact and connected. 
The relation $x\sim y$ on $M$ given by $d_M(x,y)= 0$ is symmetric 
and transitive so that $M_x=M_y$ if and only if $x\sim y$. 
So there is a well defined 
set-theoretic quotient map $\Psi : M \arrow S= M/\!\sim$. 
Note that the set $S$ is equipped with a natural metric 
induced from $d_M $. Indeed, $d_M(x',y')$ is the same for 
any points $x'\in M_x, y'\in M_y$, and hence $d_M$ induces 
a metric $d_S$ on $S$. This metric provides a topology on $S$, and since the set 
$U_{x,\epsilon} = \{ y \in M\ |\ d_M(x,y) <\epsilon \}$ is open, 
the map $\Psi: M \arrow S$ is continuous. 
Thus the metric space $S$ is also compact. 
This completes the proof of the lemma. \endproof

\hfill

\remark 
The natural quotient considered above was already proposed
in \cite{_Kobayashi:1976_} albeit little seems to be known
about its possible structure. In particular, it is known
that even when $M$ is compact, $S$ may not have the
structure of a complex variety (\cite{_Horst:2_exa_}). 
As we note in \ref{camp_10}, Campana conjectured that the 
Kobayashi metric quotient of a K\"ahler manifold has birational general type, and hence, 
a dense subset of the metric quotient should carry a complex 
(even quasi-projective) structure for such manifolds. 

\hfill

\remark 
If there is a holomorphic family of varieties 
$X_t$ smooth over a parameter space $T$ of say dimension $1$, 
then the relative construction also works by considering
the problem via that of the total space over small disks 
in $T$. In particular, there is a monodromy action on the
resulting family of compact metric spaces $S_t$ by
isometries over $T$, c.f. \S\ref{S_Ergoquotients}. 

%

\hfill

Let $M$ be a complex manifold and $h$ a hermitian metric on $M$ with its 
associated norm
$|\ |_h$.\\[-2mm]


Recall that a theorem of Royden says that the Kobayashi pseudo-metric 
$d_M$ can be obtained by 
taking the infimum of path-integrals of the infinitesimal pseudonorm 
$|\ |_M$, where 
\[ |v|_M=\inf \left \{\, \frac{1}{R}\ \ |\ \  f: D_R \rightarrow M\ 
\mbox{holomorphic}, R>0, f'(0)=v\,\right\}.
\]
Here $D_R$ is the disk of radius $R$ centred at the origin.
Recall also that $|\ |_M$ is upper-semicontinuous \cite{Siu}.

\hfill

\definition
We say that $|\ |_M$ is {\em Voisin-degenerate} at a point $x\in M$ if 
there is a sequence of holomorphic maps $\phi_n: D_{r_n}\to M$ such that 
$$\phi_n(0)\to x,\ |\phi_n^{\prime}(0)|_h=1\ \mbox{and}\ r_n\to \infty.$$

Observe that the locus $Z_M$ of $M$ consisting of points where $|\ |_M$ is 
Voisin-degenerate is a closed set. 

\hfill

\remark
If $(x,v) \in T_xM$ is a point in the tangent bundle of $M$ at $x$ which is  
Voisin-degenerate, then it does not necessarily follow that 
$|v|_M=0$, because the Kobayashi pseudometric is semicontinuous but 
might not be continuous at that point. However, the other implication 
is true: by upper semicontinuity, if $|v|_M=0$, then for any 
sequence $(x_n, v_n) \arrow (x,v)$ we have $|v_n|_M \arrow 0$, i.e., 
the point $x$ is Voisin degenerate in a strong sense. 

\hfill

The following theorem is essentially \cite[Proposition 1.19]{Vois}. 

\hfill

\theorem\label{open}
Consider the equivalence relation $x\sim y$ on $M$ 
given by $d_M(x, y)=0$ where $d_M$ is the Kobayashi
pseudo-metric on $M$. Then every non-trivial 
orbit (that is, a  non-singleton equivalence class) of this relation
consists of Voisin-degenerate points, and the union of such orbits
is a closed set. If, further, $M$ is compact, then 
each nontrivial orbit contains the image of 
a nontrivial holomorphic map $\mathbb C\to M$.

\hfill

We also need the following theorem. 

\hfill

\theorem
Assume $M$ is compact.
Then each orbit of the equivalence relation given above is connected.

\hfill

{\bf Proof:} Let $M_x$ be the orbit passing through $x$ as before and
$$M_x(n)=\left\{ \,y\in X\,\big|\, d_X(x, y)\leq \frac 1 n \,\right\}.$$ 
Then each $M_x(n)$ is compact and connected and $M_x=\cap_n M_x(n)$.
If $M_x$ is not connected, then there are disjoint open sets 
$U, V$ in $M$ separating $M_x$ leading to the contradiction
$$\emptyset = (U\cup V)^c\cap M_x=\cap_n [(U\cup V)^c\cap M_x(n)]
\neq \emptyset,$$ 
each $(U\cup V)^c\cap M_x(n)$ being nonempty compact 
as $M_x(n)$ is connected.
\endproof

\hfill

We want to exploit the existence of an automorphism
of an infinite order for the analysis of Kobayashi
metric. The following general lemma provides with
a necessary argument for a projective manifold.

\hfill

\lemma
 Let $X$ be a complex projective manifold
and $[C]$ an ample class of curves on $X$.
Let $U$ be an open domain in $X$ and $w_h$ the
volume form of a K\"ahler metric $h$
on $X$. Then for a sufficently big $n$ there is
a curve $C_1\in [nC]$ such that
$\Vol_h(C\cap U) \geq (w_h(U)/w_h(X)-\epsilon)\Vol_h(C)$
for arbitary small $\epsilon$.

\hfill

{\bf Proof:} 
The result evidently holds for $P^n$ and Fubini-Study metric
on $P^n$ since $P^n$ is homogeneous
with respect to the Fubini-Study metric.
In this case it follows from the integral volume formula
  for the family of projective lines,
parametrized by the  Grassmanian which surjects
onto $P^n$.  It immediately implies the existence
of lines which satisfy the inequality.

Similar formula holds for the family of algebraic
curves of any given degree.
In particular we obtain an infinitesimal
version of the formula which therefore
holds for any metric on projective space.
Using a finite map of an $n$-dimensional
projective manifold $X$ onto $\C P^n$
we can derive the same formula for the
K\"ahler pseudometrics induced from $\C P^n$
and then use its local nature
for any $X$.
\endproof

\hfill

\lemma \label{prev}
Let $f$ be an automorphism of infinite order on a complex 
projective manifold $X$ of dimension $n$.
Assume that there is a domain $U$ in $X$, a smooth K\"ahler metric 
$g$ on $X$ 
and positive
constants $c,c'$ such that $ cg \leq   (f^{m})^*g \leq c' g$ on $U$ for all 
powers $f^m$ of $f$. Then $f$ is an isometry of $(X, h)$ for some
K\"ahler metric $h$ on $X$ and hence some power of $f$ is contained in a
connected component of the group of complex isometries of $(X,h)$.
In particular, $X$ has a faithful holomorphic action by an abelian variety.

\hfill

{\bf Proof:} 
Let $h$ be the pull back of the Fubini-Study metric on $X$ of the embedding corresponding to a very ample line bundle $L$ on $X$.
Note that we can assume that $ag \leq   (f^{m})^*h \leq a' g$ on $\bar U$
for some  positive constants $a,a'$ which are independent of the parameter $m$.
Note that $\int_X (f^{m})^*h^n$ does not depend on $m$ since the
class of the volume $h^n$ maps into itself.
Therefore, we have
\begin{equation} \label{volume_inequality}
\mu' \int_X h^n < \int_{\bar U} (f^{m})^{*}h^n <  \mu \int_X 
h^n
\end{equation}
for some $\mu$ and $\mu'$ independent of $m$.
Let $c$ be a class of ample (i.e., very movable) curves. Then, for a sufficiently big multiple $Nc$ of the class $c$, there are curves $C\in Nc $ with $\int_{C\bigcap {\bar U}} 
 h > \nu(h,c)$, 
and similarly 
we have $\int_{C\bigcap {\bar U}}  
(f^{m})^{*}h > \nu ((f^{m})^{*}h,c)$, 
where $(h,c)$ is a pairing of the homology class $c$ and the class of kahler
metric $h$. Since
\[ a \int_{C\bigcap {\bar U}} g < \int_{C\bigcap {\bar U}} (f^{m})^{*}h < 
a' \int_{C\bigcap {\bar U}}  g,\] 
we obtain that  $((f^{m})^{*}h,c)$ is bounded from above by $a' (g,c)$
and from below by by $a (g,c)$ for  any ample
class. Since ample classes generate the dual N$_1(X)_\R$ of NS$(X)\otimes \R$
we obtain that $(f^{m})^{*}h$ as linear functional on N$_1(X)\otimes \R$
is contained in a bounded subset. 

A slightly more direct argument for
this last boundedness is as follows. Since each $f^{m*}h$ represents a K\"ahler class, it is sufficient to bound them from above as linear functionals on ample classes $c$ of curves. Note that the first inequality in \ref{volume_inequality} says that $w_{f^{m*}h}(U)> \mu' w_h(X)$ with $\mu'$ independent of $m$. By the previous lemma, one can therefore choose $C_m\in c$ such that
$$((f^{m})^{*}h,c)=\mbox{Vol}_{(f^{m})^{*}h}(C_m)\leq \delta \int_{C_m\bigcap {\bar U}} (f^{m})^{*}h < 
\delta a' \int_{C_m\bigcap {\bar U}}  g\leq \gamma (g,c),$$
where $\delta$ is independent of $m$ and $\gamma=\delta a'$.

Since there are only a finite number of integral classes in any bounded set in NS$(X)_\R$, it follows that $f^{m_0}$ leaves invariant the K\'ahler class of $h$ for some $m_0\neq 0$ and we may therefore assume that $f$ itself leaves it invariant.  
In the case of Ricci-flat $X$, $f$ must therefore be an isometry with respect to the unique Ricci-flat metric in  the K\"ahler class given by Yau's solution to the Calabi conjecture. In general,
if $H^1(X,C)=0$, then $f$ is induced from a projective action on $P^N$
under a map $X\to P^N$.  
If $H^1(X,C)\neq 0$, then we have a map from $X\times \Pic^0(X)$
to a projective family of projective spaces $\{ {\mathbb P}(H^0(L_t)^\vee), 
t\in \Pic^0(X) \}$ over $\Pic^0(X)$ and $L_t$ defines
a very ample invariant invertible sheaf on $X\times \Pic^0(X)$ over
$\Pic^0(X)$.
Hence, $f$ has to be a complex isometry on $X$ which
completes the proof of this result. 
\endproof



\hfill

\hfill

\theorem\label{_infi_auto_then_dege_Theorem_}
Let $M$ be a projective manifold with an automorphism $f$
of infinite order. Then the Kobayashi pseudo-metric $d_M$
is everywhere degenerate in the sense that $M_x\neq \{x\}$ for
all $x\in M$. Also the Kobayashi-Royden pseudo-norm $|\ |_M$
is everywhere Voisin-degenerate. Moreover,
every fiber of the map $\Psi: M\arrow S$ constructed above
contains a Brody curve and is connected.

\hfill

{\bf Proof:}
The map $f: M \arrow M$ commutes with the projection
onto $S$, and hence, induces an isometry on $S$.
Since the action of $f$ has infinite order on
$S$, there is a sequence of powers
$f^{N_i}$ which converges to the identity on $S$ by the
compactness of the group $\text{Isom}(S)$ of isometries of $S$ (in the compact 
open topology) and by
setting $N_i=n_i-n_{i-1}$ for a convergent subsequence
$f^{n_i}$ in Isom$(S)$.
We assume arguing by contradiction that $d_M$ is
non-degenerate at a point $x\in M$.
Let $U$ be the maximal subset in $M$ where $\Psi$
is a local isomorphism. Since the subsets $M_x$ are
connected, then $U$ is exactly the subset where
$\Psi$ is an embedding. The set $U$ is
invariant under $f$ and is open by \ref{open}. Hence, $f^{N_i}$ converges
to the trivial action on $U$. The boundary $\partial U$ 
of $U$ is a compact subset in $M$ with $\partial U\neq \bar U$ and
$d_M(x, \partial U) > 0$ for any point $x\in U$.
Thus, a compact subset $U_\epsilon$ which
consists of points $x\in U$ with $ d_M(x, \partial U) \geq \epsilon$
is $f$-invariant and the restriction of $d_M$ on $U_\epsilon$
is a metric. 
It is also invariant under the action of $f$ and by
theorem of Royden (\cite[Theorem 2]{_Royden:LNM_})
we know that there are smooth K\"ahler metrics $g, g'$
on $X$ with the property that $g' >  d_M \geq g $ on $U_\epsilon$.
Applying \ref{prev} we obtain that $f$
is an isometry on $M$ with respect to some K\"ahler metric.
Thus, either  $M$ has a nontrivial action of a connected 
algebraic group, and hence, trivial Kobayshi pseudometric, 
or $f$ is of finite order which contradicts our assumption. 
Thus, we obtain a conradiction also with our initial assertion 
that $d_M$ is metric on some open subset in $M$. 

Note that a limit of Brody curves is again a nontrivial
Brody curve by Brody's classical argument. By
\ref{open}, this implies that the map $\Psi: M \arrow S$ 
is everywhere degenerate, as it is degenerate in the complement 
of an everywhere dense open subset. 
\endproof

\hfill

\remark \label{camp_10}
In \cite[Conjecture 9.16]{_Campana:orbifolds_},
F. Campana conjectured that the Kobayashi quotient map of a complex projective manifold $M$
should coincide (in the birational category) with the ``core map''
of $M$, with fibers which are ``special'' and the base
which is a ``general type'' orbifold.
Then \ref{_infi_auto_then_dege_Theorem_}
would just follow, because the automorphism group of 
a general type variety is finite. Then a general fiber
of the Kobayashi quotient map contains infinitely many
points, hence its fibers are positively dimensional.

\hfill

\remark
Note that both conditions of \ref{prev} are sharp. 
It was shown by McMullen \cite{_McMullen_} that there are Kahler 
non-projective K3 surfaces with 
automoprhisms of infinite order which contain invariant
domains isomorphic to the two-dimensional ball.
There are also examples by Bedford and Kim \cite{_Bedford_Kim_} 
of rational projective surfaces $X$ with automorphisms of infinite 
order which contain an invariant ball. In this case there are 
no invariant volume forms on the variety $X$.

\section{Metric geometry of Kobayashi quotients}\label{S_Ergoquotients}


\definition
Let $M$ be a complex manifold, and $d_M$ its Kobayashi pseudometric.
Define {\bf the Kobayashi quotient} $M_K$ of $M$
as the space of all equivalence classes
$\{x\sim y \ \ |\ \ d_M(x,y)=0\}$ equipped with
the metric induced from $d_M$. 

\hfill

The main result of this section is the following theorem.

\hfill

\theorem\label{_Kob_quo_izome_Theorem_}
Let $(M,I)$ be a compact complex manifold, and $(M,J)$ its deformation.
Assume that the complex structures $I$ and $J$ are both ergodic.
Then the corresponding Kobayashi quotients are isometric.

\hfill

{\bf Proof:} Consider the limit $\lim \nu_i(I)=J$, where
$\nu_i$ is a sequence of diffeomorphisms of $M$. For each point
$x\in (M,I)$, choose a limiting point $\nu(x)\in (M,J)$
of the sequence $\nu_i(x)$. 
Fix a dense countable subset $M_0\subset M$
and replace the sequence $\nu_i$
by its subsequence in such a way that $\nu(m):=\lim \nu_i(m)$
is well defined for all $m\in M_0$.

By the upper-semicontinuity of the Kobayashi 
pseudometric, we have 
\begin{equation} \label{_decre_distances_nu^-1_Equation_}
d_{(M,J)}(\nu(x), \nu(y))\geq  d_{(M,I)}(x,y).
\end{equation}
Let $C_0$ be the union of all $\nu(x)$ for all $x\in M_0$.
Define a map $\psi:\; C_0 \arrow (M,I)$ mapping $z=\nu(x)$ to
$x$ (if there are several choices of such $x$, choose one
in arbitrary way). By \eqref{_decre_distances_nu^-1_Equation_},
the map $\psi$ is 1-Lipschitz with respect to the Kobayashi
pseudometric. We extend it to a Lipschitz map on the
closure $C$ of $C_0$. For any $x\in (M, J)$,
the Kobayashi distance between $x$ and $\psi(\nu(x))$ is
equal zero, also by \eqref{_decre_distances_nu^-1_Equation_}.
Therefore, $\psi$ defines a surjective map on Kobayashi
quotients: $\Psi:\; C_K \arrow (M,I)_K$. Exchanging $I$ and $J$, 
we obtain a 1-Lipshitz surjective map $\Phi:\; C'_K \arrow (M,J)_K$,
where $C'_K$ is a subset of $(M,I)_K$. Taking a composition
of $\Psi$ and $\Phi$, we obtain a 1-Lipschitz, surjective map
from a subset of $(M,I)_K$ to $(M,I)_K$. The following proposition
shows that such a map is always an isometry, finishing the proof of
\ref{_Kob_quo_izome_Theorem_}. \endproof

\hfill

\proposition\label{_subset_isometry_Proposition_}
Let $M$ be a compact metric space, $C\subset M$
a subset, and $f:\; C \arrow M$ a surjective 1-Lipschitz
map. Then $C=M$ and $f$ is an isometry.

\hfill

\ref{_subset_isometry_Proposition_}
is implied by the following three lemmas, 
some which are exercises found in \cite{_BBI_}.

\hfill

\lemma\label{_subse_surje_1-Lip_Lemma_}
Let $M$ be a compact metric space, $C\subset M$
a subset, and $f:\; C \arrow M$ a surjective 1-Lipschitz
map. Then $M$ is the closure of $C$.

\hfill

{\bf Proof:} Suppose that $M$ is not the closure $\bar C$ of $C$.
Take $q\in M \backslash \bar C$, and let $\epsilon = d(q, \bar C)$.
Define $p_i$ inductively, $p_0=q$, $f(p_{i+1})=p_i$.
Let $p\in \bar C$ be any limit point of the sequence $\{p_i\}$,
with $\lim_i p_{n_i}=p$. Since $f^m(p_n)\in C$ for any $m<n$,
one has $f^m(p)\in \bar C$.

Clearly, $f^{n_i}(p_{n_i})=q$.
Take $n_i$ such that $d(p, p_{n_i})<\epsilon$. 
Then $d(f^{n_i}(p), q)< \epsilon$. This is a contradiction,
because $f^n(p)\in \bar C$ and $\epsilon = d(q, \bar C)$.
\endproof

\hfill

\lemma\label{_isom_bije_Lemma_}
Let $M$ be a compact metric space, and $f:\; M \arrow M$ 
an isometric embedding. Then $f$ is bijective.

\hfill

{\bf Proof:} Follows from \ref{_subse_surje_1-Lip_Lemma_} directly.
\endproof

\hfill

\lemma\label{_1_Lip_bije_isometry_Lemma_}
Let $M$ be a compact metric space, and $f:\; M \arrow M$ 
a 1-Lipschitz, surjective map. Then $f$ is an isometry.

\hfill

{\bf Proof:} Let $d$ be the diameter of $M$, and let $K$ be the
space of all 1-Lipschitz functions $\mu:\; M \arrow [0,d]$
with the $\sup$-metric. By the Arzela-Ascoli theorem, $K$ is compact.
Now, $f^*$ defines an isometry from $K$ to itself,
$\mu \arrow \mu \circ f$. For any $z\in M$, the function
$d_z(x)=d(x,z)$ belongs to $K$. However, $d_{f(z)}$ does not
belong to the image of $f^*$ unless $d(z,x)= d(f(z), f(x))$ for all $x$,
because if  $d(z,x)< d(f(z), f(x))$, one has 
$(f^*)^{-1}(d_{f(z)})(f(x))= d(z,x)> d(f(z), f(x))$, hence
$(f^*)^{-1}(d_{f(z)})$ cannot be Lipschitz.
This is impossible by \ref{_isom_bije_Lemma_},
because an isometry from $K$ to itself must be bijective. Therefore,
the map $f:\; M \arrow M$ is an isometry.
\endproof

\hfill

The proof of \ref{_subset_isometry_Proposition_}
easily follows from \ref{_1_Lip_bije_isometry_Lemma_} and 
\ref{_subse_surje_1-Lip_Lemma_}. Indeed, by 
\ref{_subse_surje_1-Lip_Lemma_}, $f$ is a surjective, 1-Lipschitz map from
$M$ to itself, and by \ref{_1_Lip_bije_isometry_Lemma_} it is an isometry.
\endproof


\section{Eigenvalues and periodic points of hyperbolic automorphisms}\label{S_Eigen}


The following proposition follows from a simple linear-algebraic
observation.

\hfill

\proposition
Let $T$ be a holomorphic automorphism of a
hyperk\"ahler manifold $(M,I)$, and $\gamma:\; H^2(M)\arrow H^2(M)$
the corresponding isometry of $H^2(M)$. Then $\gamma$ has at most
1 eigenvalue $\alpha$ with $|\alpha|>1$, and such $\alpha$ is real.

\hfill

{\bf Proof:}
Since $T$ is holomorphic, $\gamma$ preserves the Hodge decomposition
$$H^2(M,\R)= H^{(2,0)+(0,2)}(M,\R)\oplus H^{1,1}(M,\R).$$
Since the BBF form is invariant under $\gamma$ and 
is positive definite on $H^{(2,0)+(0,2)}(M,\R)$,
the eigenvalues of $\gamma$ are $|\alpha_i|=1$ on this space.
On $H^{1,1}(M,\R)$, the BBF form has signature $(+, -,-, ..., -)$,
hence $\gamma$ can be considered as an element of $O(1,n)$.
However, it is well known that any element of $SO(1,n)$
has at most 1 eigenvalue $\alpha$ with $|\alpha|>1$, and
such $\alpha$ is real.
\endproof

\hfill

\definition
An automorphism of a hyperk\"ahler manifold $(M,I)$
or an automorphism of its cohomology algebra preserving
the Hodge type is called {\bf hyperbolic} if it acts with an eigenvalue
$\alpha$,  $|\alpha|>1$ on $H^2(M,\R)$.

\hfill

In holomorphic dynamics, there are many uses for the
{\bf $d$-th dynamical degree of an automorphism},
which is defined as follows. Given an automorphism
$T$ of a manifold $M$, we consider the corresponding
action on $H^d(M,\R)$, and  $d$-th dynamical degree is
logarithm of the maximal absolute value of its eigenvalues. In
\cite{_Oguiso:degree_}, K. Oguiso has shown that
the dynamical degree of a hyperbolic automorphism
is positive for all even $d$, and computed it explicitly
for automorphisms of Hilbert schemes of K3 which
come from automorphisms of K3. For 3-dimensional
K\"ahler manifolds, dynamical degree was computed
by F. Lo Bianco (\cite{_FLB:3-dim_}).

We compute the dynamical degree and the maximal eigenvalue
of the automorphism action on cohomology
for all even $d$ and give an upper bound
for odd ones. We also compute asymptotical growth of
the trace of the action of $T^N$ in cohomology,
which could allow one to prove that the number of
quasi-periodic points grows polynomially as the
period grows. One needs to be careful here, because
there could be periodic and fixed subvarieties,
and their contribution to the Lefschetz fixed point
formula should be calculated separately.

\hfill

\theorem\label{_eigenva_cohomo_Theorem_}
Let $(M,I)$ be a hyperk\"ahler manifold, and
$T$ a hyperbolic automorphism acting on cohomology
as $\gamma$. Denote by $\alpha$ the
eigenvalue of $\gamma$ on $H^2(M,\R)$ with $|\alpha|>1$. Then all eigenvalues of
$\gamma$ have absolute value which is a power of
$\alpha^{1/2}$. Moreover, the maximal of these eigenvalues
on even cohomology $H^{2d}(M)$ is equal to $\alpha^{d}$, and finally, 
on odd cohomology $H^{2d+1}(M)$ the maximal eigenvalue of $\gamma$ is 
strictly less than $\alpha^{\frac{2d+1}{2}}$.

\hfill

\remark Since the K\"ahler cone of $M$ is fixed by $\gamma$,
$\alpha$ is positive; see e. g.  \cite{_Canta:Milnor_}.

\hfill

\remark
From \ref{_eigenva_cohomo_Theorem_}, it follows immediately that
$\Tr(\gamma^N)$ grows asymptotically as $\alpha^{nN}$.

\hfill

We prove \ref {_eigenva_cohomo_Theorem_}
at the end of this section.

\hfill

Recall that the Hodge decomposition defines multiplicative action
of $U(1)$ on cohomology $H^*(M)$, with $t\in U(1)\subset \C$ acting
on $H^{p,q}(M)$ as $t^{p-q}$. In
\cite{_Verbitsky:coho_announce_},
the group generated by $U(1)$ for all complex
structures on a hyperk\"ahler manifold
was computed explicitly, and it was found
that it is isomorphic $G=\Spin^+(H^2(M,\R), q)$
(with center acting trivially on even-dimensional
forms and as -1 on odd-dimensional forms; see
\cite{_V:Mirror_}). Here $\Spin^+$ denotes the
connected component.

In \cite{_V:Torelli_}, it was shown that
the connected component of the group of automorphisms of $H^*(M)$
is mapped to $G$ surjectively and with compact kernel
(\cite[Theorem 3.5]{_V:Torelli_}).
Therefore, to study the eigenvalues of automorphisms
of $H^*(M)$, we may always assume that they belong
to $G$.

Now, the eigenvalues of $g\in G$ on its irreducible
representations can always be computed using the Weyl
character formula. The computation is time-consuming,
and instead of using Weyl character formula, we use
the following simple observation.

\hfill

\claim\label{_adjoint_elements_same_eigen_Claim_}
Let $G$ be a group, and $V$ its representation.
Then the eigenvalues of $g$ and $xgx^{-1}$ are equal
for all $x, g\in G$.
\endproof

\hfill

To prove \ref {_eigenva_cohomo_Theorem_}, we replace
one-parametric group containing the hyperbolic
automorphism by another one-parametric group
adjoint to it in $G$, and describe this second
one-parametric group in terms of the Hodge decomposition.

\hfill

\proposition\label{_eigenva_Proposition_}
Let $(M,I)$ be a hyperk\"ahler manifold, and
$\gamma$ an automorphism of the ring $H^*(M)$. Assume that
$\gamma$ acts on $H^{2}(M)$ with an eigenvalue $\alpha>1$.
Then all eigenvalues of
$\gamma$ have absolute value which is a power of
$\alpha^{1/2}$. Moreover, the maximal of these eigenvalues
on even cohomology $H^{2d}(M)$ is equal to $\alpha^{d}$
(with eigenspace of dimension 1), and
on odd cohomology $H^{2d+1}(M)$ 
it is strictly less than $\alpha^{\frac{2d+1}{2}}$.

\hfill

{\bf Proof:}
Denote by $G$ the group of automorphisms of $H^*(M)$.
As shown above, its Lie algebra is $\goth(so)(3,b_2(M)-3)$,
hence the connected component of $G$ is a simple Lie group.

Write the polar decomposition
$\gamma= \gamma_1 \circ \beta$, where $\gamma_1\in G$ has
eigenvalues $\alpha, \alpha^{-1}, 1, 1,..., 1$,
$\beta$ belongs to the maximal compact subgroup,
and they commute. Clearly, the eigenvalues of $\beta$
on $V$ are of absolute value 1, and absolute
values of eigenvalues of $\gamma$ and $\gamma_1$
are equal. Therefore, we can without restricting
generality assume that $\gamma=\gamma_1$ has
eigenvalues $\alpha, \alpha^{-1}, 1, 1,..., 1$.

Consider now the following one-parametric subgroup
of the complexification $G_\C\subset \Aut(H^*(M,\C))$:
$\rho(t)$ acts on $H^{p,q}$ as $t^{p-q}$, $t\in \R$.
The corresponding element of the Lie algebra has only
two non-zero real eigenvalues in adjoint action. Clearly,
all one-parametric subgroups of $G_\C=\Spin(H^2(M,\C))$
with this property are conjugate.
This implies that $\gamma$ is conjugate to an element
$\rho(\alpha)$.

By \ref{_adjoint_elements_same_eigen_Claim_},
$\gamma$ and $\rho(\alpha)$ have the same eigenvalues,
and $\rho(\alpha)$ clearly has eigenvalues
$\alpha^{\frac{d-i}2}, \alpha^{\frac{d-i-1}2}, ... \alpha^{\frac{i-d}2}$
on $H^{d}(M)$.
\endproof

\hfill

\corollary 
\[ \lim_{n\arrow \infty} \frac {\log \Tr (f^n)\restrict{H^*(M)}}n=d\log\alpha,
\]
where $2d=\dim_\C M$.
In particular, the number of $k$-periodic points grows
as $\alpha^{nk}$, assuming that they are isolated.
\endproof

\hfill

\remark
The case when $f$ admits non-isolated periodic points is treated in
\cite{_DNT:growth_}, who prove that the number of isolated $k$-periodic
points still grows no faster than $\alpha^{nk}$; the lower bound is
still unknown.

\hfill

The same argument as in 
\ref{_eigenva_Proposition_}
also proves the following theorem.

\hfill

\theorem 
Let $M$ be a hyperk\"ahler manifold, and
$\gamma\in \Aut(H^*(M))$ an automorphism 
of cohomology algebra preserving
the Hodge decomposition and acting on $H^{1,1}(M)$
hyperbolically. Denote by $\alpha$ the 
eigenvalue of $\gamma$ on $H^2(M,\R)$ with $|\alpha|>1$.
Replacing $\gamma$ by $\gamma^2$ if necessary, we may
assume that $\alpha>1$. Then all eigenvalues of
$\gamma$ have absolute value which is a power of
$\alpha^{1/2}$. Moreover, the eigenspace of
eigenvalue $\alpha^{k/2}$ on $H^d(M)$ is isomorphic 
to $H^{\frac{(d+k)}2, \frac{(d-k)}2}(M)$.
\endproof

\hfill

{\bf Acknowledgements:} Many thanks to Serge Cantat, Fr\'ed\'eric
Campana, Simone Diverio, Federico Buonerba and to 
the referee for their insightful comments 
and many references and corrections to the first version of this paper. 
The first and the fourth named authors acknowledge that the article 
was prepared within the framework of a subsidy granted to the HSE by 
the Government of the Russian Federation for the implementation of 
the Global Competitiveness Program. 

\hfill

\noindent {\sc Fedor A. Bogomolov\\
Department of Mathematics\\
Courant Institute, NYU \\
251 Mercer Street \\
New York, NY 10012, USA,} \\
\tt bogomolov@cims.nyu.edu, also: \\
{\sc National Research University, Higher School of Economics, Moscow, Russia,}
\\

\noindent {\sc Ljudmila Kamenova\\
Department of Mathematics, 3-115 \\
Stony Brook University \\
Stony Brook, NY 11794-3651, USA,} \\
\tt kamenova@math.sunysb.edu
\\

\noindent {\sc Steven Lu\\
D\'epartment de Math\'ematiques, PK-5151\\
Universit\'e du Qu\'ebec \`a Montr\'eal (UQAM)\\
C.P. 8888 Succersale Centreville H3C 3P8,\\ 
Montr\'eal, Qu\'ebec, Canada,} \\
\tt lu.steven@uqam.ca\\

\noindent {\sc Misha Verbitsky\\
{\sc Laboratory of Algebraic Geometry,\\
National Research University Higher School of Economics,\\
Faculty of Mathematics, 7 Vavilova Str. Moscow, Russia,}\\
\tt  verbit@mccme.ru}, also: \\
{\sc Universit\'e Libre de Bruxelles, D\'epartement de Math\'ematique\\
Campus de la Plaine, C.P. 218/01, Boulevard du Triomphe\\
B-1050 Brussels, Belgium\\
\tt verbit@mccme.ru
}

\end{document}